\documentclass{article}
\usepackage{amssymb,amsmath,amsthm}
\usepackage{epic}

\newtheorem{thm}{Theorem}
 \newtheorem{cor}[thm]{Corollary}
 
 \newtheorem{prop}[thm]{Proposition}
 \theoremstyle{definition}
 \newtheorem{definition}[thm]{Definition}
 \theoremstyle{remark}
 \newtheorem{rem}[thm]{Remark}
 \theoremstyle{theorem}
 \newtheorem{exmp}[thm]{Example}
\newcommand{\Aut}{\mbox{\rm Aut}}

\begin{document}

\title{Structure of $\mathbb{Z}^2$ modulo selfsimilar sublattices.\thanks{Submitted to Theoretical Computer Science}}
\author{Roberto Canogar-Mckenzie \thanks{This work was realized during a visit of this author to the
Department of Algebra and Geometry of the University of
Valladolid.} \\ {\small
 Departamento de Matem\'{a}ticas,}\\ {\small Facultad de Ciencias,
U.N.E.D.,}\\ {\small Madrid 28040, Spain } \\ {\tt
rcanogar@mat.uned.es}  \and  Edgar Mart\'{\i}nez-Moro \thanks{This
author research is partially supported by the "Junta de Castilla
y Le\'on" under Project VA-0399 and by DGICYT PB97-901.}
\\{\small Dpto. de
Matem\'{a}tica Aplicada }\\{\small Fundamental, Universidad de
Valladolid,}
\\{\small Valladolid 47002, Spain }\\ {\tt edgar.martinez@ieee.org} }
\maketitle
\begin{abstract}
In this paper we show the combinatorial structure of $\mathbb{Z}^2$
 modulo sublattices selfsimilar
to $\mathbb{Z}^2$. The tool we use for dealing with this purpose
 is the notion of  association scheme. We classify when the scheme defined by the lattice is
 imprimitive and characterize its decomposition in terms of the decomposition of the gaussian
 integer defining the lattice. This arise in the classification of different forms of tiling
 $\mathbb{Z}^2$ by lattices of this type.


\end{abstract}

\section{Introduction.}

 A {\bf similarity} $\sigma$ of {\bf norm} $c$ is a map from
 $\mathbb R^n$ to $\mathbb R^n$ such that
  $\sigma u\cdot\sigma v= c u\cdot v,\quad u,v \in \mathbb R^n$.
  Let $\Lambda$ be a two-dimensional lattice, a sublattice
  $\Lambda^{\prime} \subseteq \Lambda$ is {\bf similar} to
  $\Lambda$ if $\sigma(\Lambda) = \Lambda^{\prime}$. $\sigma$
  is also called a {\bf multiplier} of norm $c$ for $\Lambda$.
  Let us consider now the lattice $\Lambda=\mathbb Z [i]\cong \mathbb Z^2$ of Gaussian
  integers, it is a known result \cite{CRS} that the lattice $\mathbb Z^2$ have
   multipliers of norm $c$ if and only if $c=r^2+s^2,\quad r,s\in \mathbb Z$.

In this paper we shall study the combinatorial structure of
sublattices similar to $\mathbb Z [i]$ given by $(r+si)\mathbb Z [i]$, by studying
 the quotient
$\mathbb Z
[i]/(r+si)\mathbb Z [i]$.  From now on, this lattice  will be denoted $\mathbb Z
[i]_{(r+si)}$ for sort. Sublattices selfsimilar to  $\mathbb{Z}^2$ have been found useful
for recursive constructions
  of lattices (see \cite{CRS,CS}) and quotients of this lattices
  for coding two-dimensional signal constellations in coding theory
  (see \cite{Huber1,Edgar}).   We define
an {\bf association scheme} over the classes in the sublattice.
This association scheme is defined by the orbitals of a
transitive action (see \cite{Cameron} for a primer on this
constructions) . This approach is the same as defining the
association scheme given by Mannheim metric (see \cite{Edgar})
and its well known in coding theory (see   \cite[`` An all
purpouse construction'']{Sole} or \cite{Hand,Lev}). This also
arises to different
 ways of {\bf tiling} $\mathbb{Z}^2$ according the gaussian integer we
 have chosen.

The organization of the paper is as follows: in the second
section we state  some of the algebraic preliminaries underlying
the paper and also develop a general setting for dealing with
schemes derived from quotient lattices. Section 3 shows the
construction of the scheme and some of its properties such as the
expression of relation matrices. Sections 4 and 5 are devoted to
the concept of quotient schemes and its relation with tiling
 the lattice $\mathbb{Z}^2$.
 Finally in the conclusions we show some other lines of research.

\section{Algebraic Preliminaries }

\begin{definition}
Let X be a finite set. A $d$-class association scheme is a pair
$(X,(R_i)_{i \in I})$, where $I:=\{0,1, \ldots,d\}$, such that
\begin{enumerate}
\item $(R_i)_{i \in I}$ is a partition of $X \times X$,
\item $\forall i \in I, \exists j \in I$ such that $R_i^t=R_j$,
\item $R_0 := \{(x,x)|x\in X\}$,
\item there are numbers $p_{ij}^k$ such that for any pair $(x,y) \in
R_k$ the number of $z \in X$ with $(x,z) \in R_i$ and $(z,y) \in
R_j$ equals $p_{ij}^k$ and
\item $p_{ij}^k =p_{ji}^k$ for all $i,j,k \in I$.
\end{enumerate}
\end{definition}

\begin{definition} An association scheme $(X,(R_i)_{i \in I})$ is called
{\bf imprimitive} if some union of relations is an equivalence
relation distinct from the trivial ones.
\end{definition}

Let $\Lambda$ be a lattice in $n$-dimensional space $\mathbb{R}^n$
(i.e. a discrete additive subgroup of $\mathbb{R}^n$). The
\emph{automorphism group of $\Lambda$}, $\Aut(\Lambda)$, is the
set of distance-preserving transformations (or isometries) of the
space that fix the origin and take the lattice to it self. Let
$\Lambda' \subset \Lambda$ be a sublattice of $\Lambda$ and $\pi:
\Lambda \rightarrow \Lambda / \Lambda'$ be the natural projection.
We say that the function $s: \Lambda / \Lambda' \rightarrow
\Lambda$ is a \emph{section of $\pi$} if $\pi \circ s$ is the
identity.  The set $\Aut(\Lambda/ \Lambda')$ is composed of the
elements of $\Aut(\Lambda)$ that also fix $\Lambda'$.

\begin{prop}
We define an association scheme in $\Lambda / \Lambda'$ in the
following way: $(X_0, Y_0) \in \Lambda / \Lambda' \times \Lambda /
\Lambda'$ and $(X_1,Y_1) \in \Lambda / \Lambda' \times \Lambda /
\Lambda'$ are in the same relation if and only if there is a
$\sigma \in \Aut(\Lambda / \Lambda')$  such that for any section
$s: \Lambda / \Lambda' \rightarrow \Lambda$ we have that
\[\sigma(s(Y_0-X_0)) \equiv s(Y_1-X_1) (\mod \Lambda').\]
\end{prop}
\begin{proof} \begin{enumerate}
\item By the definition of the
relations, condition $1$ is satisfied.
\item The pair $(X,Y)$ is in
the same relation as $(Y,X)$ because $-(s(Y-X)) \equiv s(X-Y)
(\mod \Lambda')$. Note that the multiplication by $-1$ is an
element of $\Aut(\Lambda / \Lambda')$. It follows that $R_i^t=R_i$
for all $i \in I$. Then we say that the association scheme is
symmetric.
\item For all $X, Y \in \Lambda/ \Lambda'$, $(X,X)$ and $(Y,Y)$
are in the same relation. Thus $R_0 \supset \{(X,X)\;|\;X \in
\Lambda/ \Lambda'\}$. If $(X,Y)$ is in the same relation than
$(Z,Z)$, then $\exists \sigma \in \Aut(\Lambda/\Lambda')$ such
that:
\[\sigma(s(Y-X)) \equiv s(Z-Z) \equiv s(0) (\mod \Lambda').\]
Then
\[\sigma(s(Y-X)) \in \Lambda' \Rightarrow s(Y-X) \in
\Lambda' \Rightarrow Y-X =0\] We conclude that $R_0 =
\{(X,X)\;|\;X \in \Lambda/ \Lambda'\}$.

\item Let $(X_t,Y_t) \in R_k$ and let \[E_t=\{Z \in
\Lambda/\Lambda' \;|\; (X_t,Z)\in R_i \text{  and  } (Z,Y_t)\in
R_j\},\] for $t=0,1$. Then there is a $\sigma \in \Aut(\Lambda /
\Lambda')$  such that for any section $s: \Lambda / \Lambda'
\rightarrow \Lambda$ we have that
\[\sigma(s(Y_0-X_0)) \equiv
s(Y_1-X_1) (\mod \Lambda).\] Let $f:\Lambda \rightarrow \Lambda$
be defined as $f(x)= \sigma(x-s(X_0))+s(X_1)$. Then
$\pi(f(s(Y_0)))=Y_1$. We will prove that $\#E_0 \geq \#E_1$. For
every $Z \in E_0$ we can check that $(X_1, \pi(f(s(Z)))) \in R_i$
and also that $(\pi(f(s(Z))),Y_1)= (\pi(f(s(Z))),\pi(f(s(Y_0))))
\in R_j$. This means that $Z \in E_1$. Then $\;\pi \circ f \circ
s\;$ is a function from $E_0$ to $E_1$. We still need to show that
this function is injective. Let $Z_0, Z_1 \in \Lambda / \Lambda'$
and let us assume that
\begin{align*}
 \pi \circ f \circ s(Z_0) - \pi \circ f \circ s(Z_1)  =0
  & \Rightarrow   \pi(\sigma(s(Z_0)) - \sigma(s(Z_1)))  = 0
\Rightarrow \\
 \sigma(s(Z_0)-s(Z_1))  \in \Lambda'  & \Rightarrow s(Z_0)-s(Z_1)
    \in \Lambda'
     \Rightarrow \\  \pi(s(Z_0))-\pi(s(Z_1)) =0 & \Rightarrow
    Z_0=Z_1.
\end{align*}
\item By the symmetry of the association scheme we deduce that
      $p_{ij}^k =p_{ji}^k$.
\end{enumerate}
\end{proof}

\subsection{The Bose-Mesner algebra}

This definition in terms of relations allows us a much more
convenient way to describe association schemes in terms of
algebras.

From now on we will suppose that $X$ has $n$ ordered elements:
$X=\{x_1, \ldots, x_n\}$. If we have an association scheme
$(X,(R_i)_{i \in I})$, the family $(A_i)_{i \in I}$ of non-zero $n
\times n$ $(0,1)$-matrices will denote the adjacency matrices of
the corresponding relations (the rows and columns of $A_i$, and
all matrices of size $n \times n$ on what follows, are indexed by
$X$ in the specified order). Now we can reformulae the conditions
in terms of matrices.
\begin{enumerate}
\item $\sum_{i \in I} A_i = J$ ,
\item $\forall i \in I, \exists j \in I$ such that $A_i^*=A_j$
(where $A_i^*$ is the adjoint of $A_i$),
\item $A_0=I_n$,
\item $A_i A_j= \sum_{k \in I} p_{ij}^k A_k$ $(i,j \in I)$ and
\item $A_i A_j = A_j A_i$.
\end{enumerate}

By conditions 1 to 4, $\{A_i\}_{i \in I}$ form a base (as a vector
space over $\mathbb{C}$) of a subalgebra in $M_n(\mathbb{C})$, so
the algebra has dimension $d+1$. By 5 this subalgebra is
commutative. This algebra, $\mathcal{A}$, is called the
\emph{Bose-Mesner} algebra of the association scheme.

\section{Definition of the scheme}

Consider the lattice of Gaussian integers $\mathbb{Z}[i]$ and the
sublattice  $L=\mathbb{Z}[i]/\alpha \mathbb{Z}[i]$ the set of
Gaussian integers modulo $\alpha \mathbb{Z}[i]$ which is similar
to $\mathbb{Z}[i]$. The {\em norm} of an element $\alpha\in
\mathbb{Z}[i]$ is just $N(\alpha)=\alpha\cdot\bar\alpha$. The
{\em units} are the elements of norm $1$.  Clearly multiplication
in $L$ by an element on the group of units of the Gaussian
integers $\mathcal G=\langle i\rangle$ ($i$ is the imaginary
unit) is an isometry fixing the origin  (from now on we shall
refer to them as rotations) , and also we shall denote the group
of translations by $\mathcal T$.

In the following discussion we will need some notation on
permutation groups acting on finite sets. For an reference on
this topic see \cite{Alperin,Cameron}. For a given  permutation
group $\mathcal G$ of elements of a finite set $X$ we denote the
orbit of an element $x\in X$ as $({\mathcal G})(x)=\{ xg\mid g\in
{\mathcal G}\}$. Two orbits are either identical or disjoint.  We
denote by ${\mathcal O}({\mathcal G}\mid X)$ the set of all the
orbits of the action. We denote the stabilizers  by ${\mathcal
G}_x =\{g\in {\mathcal G}\mid x=xg \}$. It is well known the
relation between orbits and stabilizers given by:
\[
\begin{array}{ccc}
\Bigl({\mathcal G}\Bigr)(x) & \to & {\mathcal G}/{\mathcal G}_x\\
xg &\mapsto &{\mathcal G}_x g
\end{array}
\]
The action is {\em transitive} if there is just one orbit. If the action is transitive,  a {\em congruence} is a $\mathcal G$-invariant equivalence relation on $X$. We say the action is {\em imprimitive} if has a non-trivial equivalence.

Consider now the {\em semi direct} product of both groups
${\mathcal{H}}=\mathcal{G} \ltimes {\mathcal{T}}$. Roughly
speaking we will also denote by ${\mathcal{H}}$ the permutation
group on $L$ generated by the permutation $(\alpha\mapsto
i\alpha)$ and the translations in ${\mathcal{T}}$. It is clear
that ${\mathcal{H}}$ acts transitively on $L$.

 Consider the orbits of the action:
 $$ \mathcal H\times (L\times L)\to (L\times L)$$
induced by the action of ${\mathcal{H}}$ on $L$. They are called
{\em orbitals} and they are the relations of an association scheme
\cite{Cameron}. Since $\mathcal{H}$ is a transitive permutation
group, if we take $G_{0}$ (i.e. those permutations that fix $0$),
it is well known  that we have the coset decomposition: $
{\mathcal{H}}=\bigl( G_0\bigr) (p_{0})\cup \dots\cup
\bigl(G_0\bigr)(p_t) $, where $p_i$ is the permutation
transforming $0$ in some complex number $\beta$ belonging to the
coset. Therefore  orbits can be rewritten as:
\begin{equation}\label{G0}
(x,y)\in R_k \Leftrightarrow x-y\in \bigl(G_0\bigr)( p_k)
\end{equation}
In our case $G_0=\mathcal G$.

\begin{exmp}\label{ex:1} Consider the sublattice $\mathbb Z
[i]_{(2+2i)}$. We can give a system of representants in the
fundamental region given in the figure 1.

\begin{figure}[h]\label{fig:1}\caption{$\mathbb Z
[i]_{(2+2i)}$}
\begin{center}
\setlength{\unitlength}{0.65mm}
\begin{picture}(50,50)
\put(0,20){\vector(1,1){20}} \put(22.5,42.5){\makebox(0,0){{\tiny
2+2i}}} \put(0,20){\vector(1,-1){20}} {\dottedline(-5,20)(55,20)}
{\dottedline(0,0)(0,40)} {\dottedline(20,40)(40,20)(20,0)}
\multiput(0,20)(10,0){4}{\circle*{2}}\multiput(40,20)(10,0){2}{\circle{2}}
\multiput(10,30)(10,0){2}{\circle*{2}}\multiput(30,30)(10,0){3}{\circle{2}}
\multiput(10,10)(10,0){2}{\circle*{2}}\multiput(30,10)(10,0){3}{\circle{2}}
\multiput(0,0)(10,0){6}{\circle{2}}
\multiput(0,40)(10,0){6}{\circle{2}}
\put(0,10){\circle{2}}\put(0,30){\circle{2}}
\end{picture}
\end{center}
\end{figure}
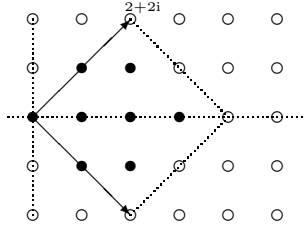

 The orbitals of the previous action are given by the cosets (see (\ref{G0})):
 \begin{align*}
 \bigl(G_0\bigr)\cdot
(0)= \{ 0\}\qquad & \bigl(G_0\bigr)\cdot (1)= \{ \pm 1, \pm i\} \\
\bigl(G_0\bigr)\cdot (2)= \{ 2\}\qquad & \bigl(G_0\bigr)\cdot
(1+i) = \{ 1+i, 1-i\} \end{align*} Which allows us construct the
relations in the association scheme.
\end{exmp}

\begin{thm}\label{th:1} The association scheme $(X,(R_i)_{i \in I})$ with $X=L$ and
$(R_i)_{i \in I}$ defined by the orbitals of the action above is
primitive if and only if $\alpha$ is a prime in $\mathbb{Z}[i]$.
\end{thm}
\begin{proof} Suppose that $\alpha$ is not a prime in
$\mathbb{Z}[i]$, i.e. $\alpha =\alpha_1\cdot\alpha_2$. Then it is
clear that the equivalence in $L$ given by the quotient
$L/[\mathbb{Z}[i]/ \alpha_1\mathbb{Z}[i]]\subset L$ is
$\mathcal{G}$-invariant, and therefore the action is imprimitive,
hence by proposition 2.9.3 in \cite{BCN} the association scheme is
imprimitive.

Conversely, if $\alpha$ is prime then is a well known fact in number theory (see \cite{Hardy})
 that either it is $1+i$ multiplied by an unit or:
\begin{enumerate}
\item $N(\alpha )=p \equiv 1\mod 4$, $p$ an odd prime, and in this case the lattice $L$ has $p$ points, and clearly $\mathcal{T} \cong \mathbb{Z}_p$.
\item $\alpha = p\in \mathbb{Z}$ and $p\equiv 3 \mod 4$, in this case $L$ can be represented as $\mathbb{Z}_p[i]$, and therefore has $p^2$ points and $\mathcal{T}\cong \mathbb{Z}_p\times\mathbb{Z}_p$.
\end{enumerate}
In the  case $1+i$ it is obvious that the scheme is primitive, in
the other two cases we have that the estabiliser of any point
$\beta$ in the lattice is given by the group of rotations around
the point (i.e. $G_\beta=t_\beta\mathcal{G} (t_\beta)^{-1}$, where
$t_\beta$ is the translation of vector $\beta$ ). It is clear
that in both cases above there is no group between
$\mathcal{G_\beta}<\mathcal{H}$ since the group generated by a
single element in $\mathcal T$ and the group $\mathcal G_\beta$
is $\mathcal H$ (See remark below for a further explanation).
Therefore the estabiliser is maximal, and by theorem 1.9 in
\cite{Cameron}, the action is primitive, and so is the
association scheme.
\end{proof}

\begin{rem} Indeed, if $G_\beta$ are the rotations around the point $\beta$ and we add a new rotation around other point, say $\beta^\prime\neq\beta$, then the group generated contains all rotations around $\beta^\prime$. Moreover, the composition of the $i$  rotation  arund one of the points, and $i^3$ around the other is a translation. So in the proof above we can suppose that we always add one translation.

When there is  a prime number of points, it is clear that a single
translation generates the group $\mathcal T$ since is cyclic of
prime order. In the case $\mathcal T\cong
\mathbb{Z}_p\times\mathbb{Z}_p$ a translation $t$ and the
translation  $i^{-1}\circ t\circ i$ are independent and  generate
the whole group $\mathcal T$. In remark \ref{rem:codes} we will
show the relation of this facts with the construction of
constellations representing $GF(p)$ and $GF(p^2)$ respectively.

\end{rem}
 The association scheme defined is a {\bf {trans\-la\-tion} {aso\-cia\-tion} {sche\-me}}, and
 in the case of a prime number of points (i.e. $|X|$ is a prime) is a {\bf cyclotomic} scheme
 (see \cite{BCN} for the definitions).
\begin{rem}\label{rem:codes}
The two constructions we present below  are used to construct
two-dimensional modulo metrics (in particular Mannheim distance)
for coding theory (for further details see
\cite{Huber1,Huber2,Edgar}).
\begin{itemize}
\item[a)] Let
$\pi\in\mathbb{Z}[i]$ be an element whose norm is a prime integer $p$, and $ p\equiv 1 \mod 4$. It is well known (Fermat's two
square theorem) that $p$ can be written as:
\[
p=a^2 +b^2 =\pi\bar\pi\hbox{ where } \pi=a+ib \hbox{ (not
unique)}.
\]
If we denote by $\mathbb{Z}[i]_\pi$ the set of Gaussian integers
modulo $\pi$, we define the modulo function $\nu
:\mathbb{Z}[i]\to\mathbb{Z}[i]_\pi$  associating to each class in
$\mathbb{Z}[i]_\pi$ its representant with smallest norm :
\begin{equation}\label{rounding}
\nu (\xi )=r \hbox{ where } \xi =q \pi +r\hbox{ and } \|
r\|=\min\{\|\beta\| \mid \beta =\xi  \rm{mod} \pi \}
\end{equation}
This can be done because $\mathbb{Z}[i]$ is and Euclidean domain.
The quotient $q$ can be calculated as $[\frac{\alpha\bar\pi}{p}]$
where $[x ]$ denotes the Gaussian integer with closest real and
imaginary part to $x$. The quotient $q$ can be calculated as
$[\frac{\alpha\bar\pi}{p}]$ where $[x ]$ denotes the Gaussian
integer with closest real and imaginary part to $x$.

 Taking  the carrier set of $GF(p)$ as $\{0,1,\dots ,p-1\}\subset \mathbb{Z}$,
 we can restrict to $GF(p)$ the application $\nu$  so that it induces
an isomorphism  $\nu :GF(p)\to \mathbb{Z}[i]_\pi$ given by:
\[
\hbox{For }g\in GF(p)\quad\nu
(g)=g-\bigl[\frac{g\bar\pi}{p}\bigr]\pi
\]
So $GF(p)$ and $\mathbb[i]_\pi$ are mathematically equivalent and
we shall use from now on that carrier set for a sort notation.
\item[b)] In the case {$p\equiv 3  \mod   4$}
$\pi =p\in \mathbb Z$ and the isomorphism above does not apport any relevant information over $GF(p)$.
For this type of primes $-1$ is a quadratic non residue of $p$, hence we get the following isomorphism between $GF(p^2)$ and ${\mathbb{Z}_p[i]}$ where :
\[
{\mathbb{Z}_p[i]} =\Bigl\{ k+il\mid k ,l\in \bigl\{\frac{-(p-1)}{2},\dots ,-1,0,1,\dots , \frac{(p-1)}{2}\bigr\}\Bigr\}
\]
\end{itemize}
\end{rem}
\begin{exmp}
Consider ${\mathbb{Z}[i]}_{3+2i}$ and ${\mathbb{Z}_7[i]}$, given
by the carrier sets defined as in the previous remark . We have
an alternative pictorial representation of them to usual one
derived as in  example \ref{ex:1} given by figure 2. This
representation is more suitable for showing the symmetries and
rotations within the fundamental region. For the association
scheme of this constellations of points see example \ref{ex:cod}.
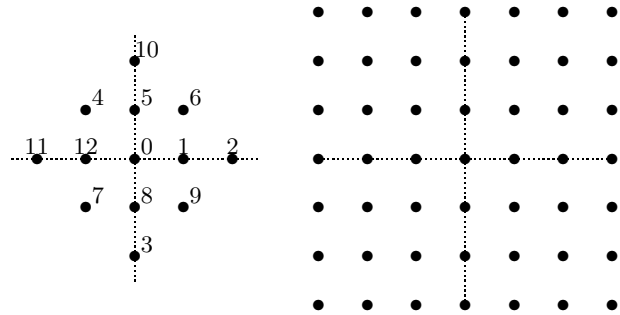
\begin{figure}[h]\label{fig:codes}\caption{$GF(13)$ as ${\mathbb{Z}[i]}_{3+2i}$  and $GF(7^2)$ as ${\mathbb{Z}_7[i]}$}
\begin{center}
\setlength{\unitlength}{0.65mm}
\begin{picture}(55,70)
\multiput(10,30)(10,0){5}{\circle*{2}}
\multiput(20,20)(10,0){3}{\circle*{2}}
\multiput(20,40)(10,0){3}{\circle*{2}} \put(30,10){\circle*{2}}
\put(30,50){\circle*{2}}
{\dottedline(5,30)(55,30)}
{\dottedline(30,5)(30,55)} \put(10,32.5){\makebox(0,0){{\small
11}}} \put(20,32.5){\makebox(0,0){{\small 12}}}
\put(32.5,32.5){\makebox(0,0){{\small 0}}}
\put(40,32.5){\makebox(0,0){{\small 1}}}
\put(50,32.5){\makebox(0,0){{\small 2}}}
\put(22.5,42.5){\makebox(0,0){{\small 4}}}
\put(32.5,42.5){\makebox(0,0){{\small 5}}}
\put(42.5,42.5){\makebox(0,0){{\small 6}}}
\put(32.5,52.5){\makebox(0,0){{\small 10}}}
\put(22.5,22.5){\makebox(0,0){{\small 7}}}
\put(32.5,22.5){\makebox(0,0){{\small 8}}}
\put(42.5,22.5){\makebox(0,0){{\small 9}}}
\put(32.5,12.5){\makebox(0,0){{\small 3}}}
\end{picture}$\qquad$
\begin{picture}(60,60)
{\dottedline(0,30)(60,30)} {\dottedline(30,0)(30,60)}
\multiput(0,30)(10,0){7}{\circle*{2}}
\multiput(0,20)(10,0){7}{\circle*{2}}
\multiput(0,10)(10,0){7}{\circle*{2}}
\multiput(0,0)(10,0){7}{\circle*{2}}
\multiput(0,40)(10,0){7}{\circle*{2}}
\multiput(0,50)(10,0){7}{\circle*{2}}
\multiput(0,60)(10,0){7}{\circle*{2}}

\end{picture}
\end{center}
\end{figure}

\end{exmp}

\subsection{Matrix expression of the algebra}

In this section we develop and easy way for describing the
matrices of the Bose-Mensner algebra associated to this
sublattices in terms of circulant matrices.

\begin{definition} Let $M$ be a $n\times n$ matrix and let $\{ a_i\}_{i=0,\dots ,n-1}$
be the first row of $M$.  Then $M$ is {\bf circulant} if:
\[M_{ij}=a_{(i-j)\mod n},\quad i,j=0,1,\dots, n\]
\end{definition}
First we shall insight in the case of a prime number of points
$p$. In this case we have that the translations $\mathcal T$ are
a clyclic group of prime order, hence, if we choose an element
generating $\langle t\rangle=\mathcal T$, any element $l\in L$
can be rewritten as $l=t^j(0),\quad 0\leq j\leq p-1$, and
choosing the order given by $j$ for the points in the scheme it is
clear that the matrices $R_i$ are circulant since $t$ is an
isometry.

\begin{exmp}\label{ex:cod} In the case $L={\mathbb{Z}[i]}_{3+2i}$ of the previous
example the orbits can be described knowing the cosets:
\begin{align*}\bigl(G_0\bigr)\cdot (p_0)= \{ 0\} \qquad & \bigl(G_0\bigr)\cdot (p_1)=
\{ 1,5,12,8\} \\ \bigl(G_0\bigr)\cdot (p_2)= \{ 6,4,7,9\} \qquad &
\bigl(G_0\bigr)\cdot (p_3) = \{ 2,10,11,3\}\end{align*} And if we
choose the translation given by adding $1$ to each point for
ordering the relations (the usual order for $GF(p)$), then they
are given by:
\begin{align*}
D_0 &=[1,0,0,0,0,0,0,0,0,0,0,0,0]\quad
D_1 =[0,1,0,0,0,1,0,0,1,0,0,0,1]\\
D_2 &=[0,0,0,0,1,0,1,1,0,1,0,0,0]\quad D_3
=[0,0,1,1,0,0,0,0,0,0,1,1,0]
\end{align*}
Note that each matrix $D$ is represented only by its first row
since they are circulant. Moreover, as usual we shall collect al
the information in a single matrix as: $$
[0,1,3,3,2,1,2,2,1,2,3,3,1]$$
\end{exmp}
\begin{rem} The eigenvalues of a circulant matrix are computed
easily (see \cite{Biggs}) as sums of roots of the unit, and
hence, so are the eigenvalues of the scheme (see \cite{Edgar}).
\end{rem}

In the case of a non-prime number of points, we have a slightly
modified construction based on the decomposition of $\mathcal T$
as the direct product of cyclic subgroups. The idea is the same
as in the previous example, and now the matrices are circulant in
blocks given by the clyclic subgroups. We shall illustrate this
idea with an example.
\begin{exmp}\label{ex:1cont} Consider $L=\mathbb Z
[i]_{(2+2i)}$ as in example \ref{ex:1}. Consider the isomorphism
given by:
\begin{align*}
L=\mathbb Z [i]_{(2+2i)} & \stackrel{\sim}{\longrightarrow}
\mathbb{Z}_4\times \mathbb{Z}_2\\
1 &\mapsto (1,0)\\
1+i &\mapsto (0,1)
\end{align*}
Consider now the elements in the order given by
$$(0,0),(1,0),\dots,(3,0),(0,1),(1,1),\dots, (3,1)$$ The relation
matrices can be expressed now as block circulant matrices: $$
D_i=\left[ \begin{array}{cc} D_{i1} & D_{i2}\\ D_{i2} &
D_{i1}\end{array}\right]\quad i=0,\dots ,3$$ Where each block is
a circulant matrix. The relations in example \ref{ex:1} can be
represented in the same fashion as previous example: $$[0,1,2,1 |
3,1,3,1]$$ where $|$ means the division of the blocks.
\end{exmp}
\begin{prop} The relation matrices can be expresed as
block-circulant matrices, where each block is also a circulant
matrix.
\end{prop}
\begin{proof}
It follows directly  from  the reasoning above, i.e. the decomposition  of the group of
translations in a product of cyclic groups.
\end{proof}

\section{Quotient schemes}

We introduce in this section the concept of quotient scheme. For
an account on this topic see \cite{BCN}. Suppose a set of indices
$\widetilde{0}$, and let us define the equivalence relation $\sim$
among the set of indices of the relations in the scheme
$({X},{\mathcal R})$ as follows: $$ a\sim b \hbox{ if }
p_{ab}^i\neq 0 \hbox{ for some } i\in \widetilde{0}$$ As usual
$\widetilde{0}$ is the class of $0$ and we write $\widetilde{a}$
for the relation containing $a$.

\begin{definition} We define a {\bf quotient scheme} $(\tilde{X},
\tilde{\mathcal R})$ of $({X},{\mathcal R})$ with respect to
$\tilde{0}$ as the association scheme whose point set is the set
$\tilde{X}$ of equivalence classes on $X$, and whose relations
are $\tilde{R}_{\tilde{i}}$, with
\[\tilde{R}_{\tilde{i}}=\{(\tilde{x},\tilde{y})
\mid \hbox{ \rm for } x\in \tilde{x},y\in \tilde{y}\hbox{ \rm we
have } (x,y)\in R_i\hbox{ \rm with } i\in \tilde{\i}\}\]
\end{definition}
\begin{prop} If the scheme defined on $L$ is imprimitive we can
define a quotient scheme where $\tilde{0}$ is given by  the
classes of some of the elements  divisors of $0$.
\end{prop}
\begin{proof} It is obvious by the first part of the proof of theorem in \ref{th:1}.
\end{proof}
\begin{exmp} Consider $L=\mathbb Z
[i]_{(2+2i)}$ as in examples \ref{ex:1} and \ref{ex:1cont}, and
consider the relation given by $\tilde{0}=\{ 0,2\}$. With the
notation in example \ref{ex:1cont} we have:
\begin{align*}
\widetilde{(0,0)}=\{(0,0),(2,0)\} &\quad
\widetilde{(1,0)}=\{(1,0),(3,0)\}\\
\widetilde{(0,1)}=\{(0,1),(2,1)\} &\quad
\widetilde{(1,1)}=\{(1,1),(3,1)\}\end{align*} And the relation
matrix in the order given by
$\widetilde{(0,0)},\widetilde{(1,0)},\widetilde{(0,1)},\widetilde{(1,1)}$
is:
\[ [0,1,2,1]\]
\end{exmp}
\begin{rem} Indeed, in the previous example, the translations are given by the
group $\mathbb{Z}_2\times \mathbb{Z}_2$, (this
can be seen also because $2=(1+i)(1-i)$ is not a gaussian prime)
and admits a further quotient scheme given by:
\[ \widetilde{\widetilde{(0,0)}}=\{\widetilde{(0,0)},\widetilde{(0,1)}\} \quad
\widetilde{\widetilde{(1,0)}}=\{\widetilde{(1,0)},\widetilde{(1,1)}\}\]
This corresponds with the identifications in the figure 3.
\begin{figure}[h]\caption{Quotients of $\mathbb Z
[i]_{(2+2i)}$}
\begin{center}
\setlength{\unitlength}{0.65mm}
\begin{picture}(50,50)
\put(0,20){\vector(1,1){20}} \put(22.5,42.5){\makebox(0,0){{\tiny
2+2i}}} \put(0,20){\vector(1,-1){20}} {\dottedline(-5,20)(55,20)}
{\dottedline(0,0)(0,40)} {\dottedline(20,40)(40,20)(20,0)}
\multiput(0,20)(10,0){4}{\circle*{2}} \put(20,30){\circle{2}}
\put(20,10){\circle{2}} \put(10,30){\circle*{2}}
\put(10,10){\circle*{2}}
\put(0,20){\circle{3}}\put(20,20){\circle{3}}
\put(10,20){\circle{3}}\put(30,20){\circle{3}}
\put(10,20){\circle{5}}\put(30,20){\circle{5}}
{\dottedline(5,25)(25,25)(25,5)(5,5)(5,25)}
\end{picture}\qquad
\begin{picture}(50,50)
\put(0,20){\vector(1,1){20}} \put(22.5,42.5){\makebox(0,0){{\tiny
2+2i}}} \put(0,20){\vector(1,-1){20}} {\dottedline(-5,20)(55,20)}
{\dottedline(0,0)(0,40)} {\dottedline(20,40)(40,20)(20,0)}
\multiput(0,20)(10,0){4}{\circle*{2}}
\multiput(10,30)(10,0){2}{\circle*{2}}
\multiput(10,10)(10,0){2}{\circle*{2}} \put(20,30){\circle{3}}
\put(30,20){\circle{3}} \put(10,20){\circle{3}}
\put(20,10){\circle{3}} {\dottedline(0,15)(30,15)}
{\dottedline(5,25)(25,25)(25,5)(5,5)(5,25)}
\end{picture}
\end{center}
\end{figure}
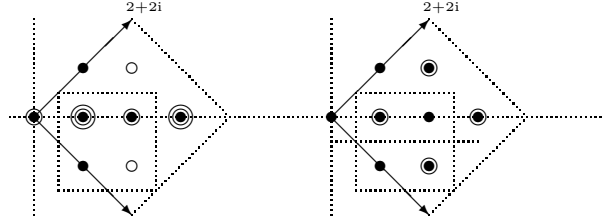
\end{rem}
\begin{rem}\label{rem:involution} Note that following \cite[pp.52]{BCN}, if we let $A=\{a\mid
p_{aa}^0=1\}$ then for each $a\in A\setminus \{ 0\}$ we find a
involution $\sigma_a:x\mapsto \bar x,\quad (x,\bar x)\in R_a$,
and if we set $\sigma_0=1$ them
$\sigma_a\sigma_b=\sigma_b\sigma_a= \sigma_c$ where $c$ is
determined by $p_{ab}^c\neq 0$. $A$ clearly has the structure of
an elementary abelian 2-group.
\end{rem}
\begin{rem} There is a sort of Jordan-H\"older theory relating the
facts above that for our example we can summarize as follows:

\begin{center} $\begin{array}{c}

\mathbb{Z}_{(2+2i)}\\
\bigcup\\
\mathbb{Z}_{(2)}\\
\bigcup\\
\mathbb{Z}_{(1+i)}\end{array}\qquad$
\begin{tabular}{|c|c|c|}
  \hline
  $\mathbb{Z}_{(\alpha)}$ & $\mathcal T$ & $\{a\mid p_{aa}^0=1\}$ \\ \hline\hline
  $\mathbb{Z}_{(2+2i)}$ & $\mathbb{Z}_2\times\mathbb{Z}_2\times\mathbb{Z}_2$ & $0,2$  \\ \hline
  $\mathbb{Z}_{(2)}$ & $\mathbb{Z}_2\times\mathbb{Z}_2$ & $\widetilde{(0,0)},\widetilde{(0,1)}$  \\ \hline
  $\mathbb{Z}_{(1+i)}$ &  $\mathbb{Z}_2$ & $\widetilde{\widetilde{(0,0)}}$
  \\
   \hline
\end{tabular}
\end{center}
\end{rem}

\section{Relation with tilings}

Up to now we have shown that we have tree types of primitive
sublattices, depending on the gaussian prime defining it, this
arises to three different forms of tiling the whole
$\mathbb{Z}^2$.
\begin{itemize}
\item With tiles of type $\mathbb{Z}[i]_{(1+i)}$;
\item with tiles of type $\mathbb{Z}[i]_{(a+bi)}$ where
$N(a+bi)=p$ an odd prime $\equiv 1 \mod 4$;
\item with tiles of type $\mathbb{Z}_p[i]$; $p$ an odd prime $\equiv 3 \mod
{4}$.
\end{itemize}

Indeed this has a close relationship with the well known fact in
number theory given by the criterion for representability of a
number $N$ by the sum of  two squares which says that any prime
factor of $N$ which is of the form $4k+3$ must divide $N$ to an
even power exactly \cite{Davenport}. In our setting this means
(as we have seen) that the regions defined by an integer of norm
$p^n$, $p$ a prime of this type must have $p^{2n}$ points. This
fact relates each prime in the factorization of the norm $N$ with
a type of the primitive tiles above.

In the three cases above the boundary of the Voronoi cell of the
sublattice $\alpha \mathbb{Z}[i]$ does not contain any extra
points of $\mathbb{Z}^2$ so following notation in \cite{BCN} we
called such sublattices {\em clean}. Hence, in the case of clean
sublattices we can find a complete set of representatives of the
non-zero classes within the region bounded by $\alpha$ and
$i\cdot \alpha$.  Moreover they are clean if there is an odd
number of points \cite{BCN}, i.e. there is no an involution (see
remark \ref{rem:involution}), so it follows directly the
following result (as expected):
\begin{cor}
Any quotient of a clean sublatice $L$ defined as in section 3 is
also clean.
\end{cor}

\begin{rem} The primitive schemes above are the finest translation schemes from our
setting. Recall that in the schemes defined for an odd prime we
have in all cases all orbits are of size $4$, but the schemes are
not pseudocyclic (see \cite{BCN} for a definition) since $\sum_i
p_{ii}^k\neq 3$. We can go a bit farther with the following result
of Rao, Ray-Chau and Singhi \cite[pp.52 ii]{BCN},\cite{RRS}: they
state that the finest translation association scheme for a set of
odd order is pseudocyclic and the other translation schemes for
the same set can be derived  from this one by merging classes. In
the case  of the schemes in this paper the scheme they recall is
the one generated by the $1,i^2$ rotations an the translations,
and clearly each relation of our scheme arises from the merging of
two relations of it.
\end{rem}

\section{Conclusions}

We have study the combinatorial structure of the association
schemes derived from sublattices given by $\mathbb Z
[i]/(r+si)\mathbb Z [i]$. As we have seen, there are close
connections of this type of lattices with coding theory, recursive
construction of lattices and self-similar lattices
\cite{BCN,Huber1,Huber2,Edgar}. In the study of the scheme
defined plays a central role the factorization of the gaussian
integer $(r+si)$ and also the factorization of the order of the
group of translations $\mathcal T$ (i.e. the number of points in
the sublattice). We have characterize the primitive cases and also
identify the cases where the Voronoi cell is clean both from a
combinatorial point of view and from a number theoretical one. We
can see the primitive case as the finest tiles of the lattice
$\mathbb Z^2$ and they are useful in coding two dimensional
signal spaces with Mannheim metric \cite{Huber1}. We have also
shown how derive an easy expresion of the matrices defining the
relations in the schemed based in thier circulant structure. A
similar study can be done with hexagonal schemes an hexagonal
metric based on Einsestein-Jacobi integers \cite{Edgar} and will
be shown elsewhere. Future trends of investigation point toward
classifying the partition designs derived from this type of
association schemes \cite{Edgar2}.

\end{document}